\documentclass[11pt]{amsart2000}
\usepackage{amssymb}
\newtheorem{theorem}{Theorem}[section]
\newtheorem{corollary}[theorem]{Corollary}
\newtheorem{lemma}[theorem]{Lemma}
\theoremstyle{definition}
\newtheorem{question}[theorem]{Question}
\newtheorem{ack}[theorem]{Acknowledgements}
\newtheorem{discussion}[theorem]{Discussion}
\newtheorem{remarks}[theorem]{Remarks}
\newtheorem{remark}[theorem]{Remark}
\newtheorem{example}[theorem]{Example}
\newtheorem{notation}[theorem]{Notation}
\newtheorem{definition}[theorem]{Definition}
\newtheorem*{mydefinition}{Definition}
\def\cc {{\mathfrak c}}
\def\mm {{\mathfrak m}}

\def\PP {{\mathbb P}}
\def\QQ {{\mathbb Q}}
\def\RR {{\mathbb R}}
\def\TT {{\mathbb T}}
\def\ZZ {{\mathbb Z}}
\def\sB {\mathcal{B}}
\def\sI {\mathcal{I}}

\def\sK {\mathcal{K}}
\def\sN {\mathcal{N}}
\def\sP {\mathcal{P}}
\def\sT {\mathcal{T}}
\def\cf {\mathrm{cf}}

\def\Hom {\mathrm{Hom}}
\def\non {\mathrm{non}}

\def\tor {\mathrm{tor}}

\begin{document}
\setlength{\unitlength}{0.01in}
\linethickness{0.01in}
\begin{center}
\begin{picture}(474,66)(0,0)
\multiput(0,66)(1,0){40}{\line(0,-1){24}}
\multiput(43,65)(1,-1){24}{\line(0,-1){40}}
\multiput(1,39)(1,-1){40}{\line(1,0){24}}
\multiput(70,2)(1,1){24}{\line(0,1){40}}
\multiput(72,0)(1,1){24}{\line(1,0){40}}
\multiput(97,66)(1,0){40}{\line(0,-1){40}}
\put(143,66){\makebox(0,0)[tl]{\footnotesize Proceedings of the Ninth Prague Topological Symposium}}
\put(143,50){\makebox(0,0)[tl]{\footnotesize Contributed papers from the symposium held in}}
\put(143,34){\makebox(0,0)[tl]{\footnotesize Prague, Czech Republic, August 19--25, 2001}}
\end{picture}
\end{center}
\vspace{0.25in}
\setcounter{page}{23}
\title{Concerning the dual group of a dense subgroup}
\author{W. W. Comfort}
\address{W. W. Comfort\\
Department of Mathematics, Wesleyan University\\ 
Middletown, CT 06459}
\thanks{The first author was an invited speaker at the Ninth Prague 
Topological Symposium.}
\email{wcomfort@wesleyan.edu}
\author{S. U. Raczkowski}
\author{F. Javier Trigos-Arrieta}
\address{S. U. Raczkowski and F. Javier Trigos-Arrieta\\
Department of Mathematics, California State University, Bakersfield\\
9011 Stockdale Highway\\
Bakersfield, CA, 93311-1099}
\email{racz@cs.csubak.edu}
\email{jtrigos@cs.csubak.edu}
\begin{abstract}
Throughout this Abstract, $G$ is a topological Abelian group and  
$\widehat{G}$ is the space of continuous homomorphisms from $G$ into 
$\mathbb T$ in the compact-open topology. A dense subgroup $D$ of $G$ {\it 
determines} $G$ if the (necessarily continuous) surjective isomorphism 
$\widehat{G}\twoheadrightarrow\widehat{D}$ given by $h\mapsto h|D$ is a 
homeomorphism, and $G$ is {\it determined} if each dense subgroup of $G$ 
determines $G$. The principal result in this area, obtained independently 
by {\sc L. Au}{\ss}{\sc enhofer} and {\sc M.~J.~Chasco}, is the following: 
Every metrizable group is determined. 
The authors offer several related results, including these.
\begin{enumerate}
\item There are (many) nonmetrizable, noncompact, determined groups.
\item If the dense subgroup $D_i$ determines $G_i$ with $G_i$ compact, 
then $\oplus_i\,D_i$ determines $\Pi_i\,G_i$. In particular, if each $G_i$ 
is compact then $\oplus_i\,G_i$ determines $\Pi_i\,G_i$.
\item Let $G$ be a locally bounded group and let $G^+$ denote $G$ with its Bohr 
topology. Then $G$ is determined if and only if ${G^+}$ is determined. 
\item Let $\non(\mathcal{N})$ be the least cardinal $\kappa$ such that some 
$X \subseteq \TT$ of cardinality $\kappa$ has positive outer measure. No 
compact $G$ with $w(G)\geq\non(\mathcal{N})$ is determined; thus if 
$\non(\mathcal{N})=\aleph_1$ (in particular if CH holds), an infinite 
compact group $G$ is determined if and only if $w(G)=\omega$.
\end{enumerate}
Question. 
Is there in ZFC a cardinal $\kappa$ such that a compact group $G$ is  
determined if and only if $w(G)<\kappa$? Is $\kappa=\non(\mathcal{N})$?
$\kappa=\aleph_1$?
\end{abstract}
\thanks{Parts of this paper appeared in~\cite{raczphd}.
Furthermore, portions of this paper were presented by the authors at the Ninth 
Prague Topological Symposium (Praha, August, 2001) and at the 2002 Annual
Meeting of the American Mathematical Society (San Diego, January, 2002).
A full treatment, with proofs, will appear elsewhere~\cite{crt}.}
\thanks{W. W. Comfort, S. U. Raczkowski, and F. Javier Trigos-Arrieta,
{\em Concerning the dual group of a dense subgroup},
Proceedings of the Ninth Prague Topological Symposium, (Prague, 2001),
pp.~23--35, Topology Atlas, Toronto, 2002; {\tt arXiv:math.GN/0204147}}
\maketitle

\setcounter{section}{-1}
\section{Terminology, Notation and Preliminaries}

For $X$ a set and $\kappa$ a cardinal, we write
$[X]^\kappa=\{A\subseteq X:|A|=\kappa\}$.

For each space $X=(X,\sT)$ we write 
$$\sK(X):=\{K\subseteq X:K \mbox{ is $\sT$-compact}\}.$$
All groups considered here, whether or not equipped with a
topology, are Abelian groups written additively. The identity of a group
$G$ is denoted $0$ or $0_G$, and the torsion subgroup of $G$ is denoted
$\tor(G)$. The reals, rationals, and integers are
denoted $\RR$, $\QQ$, and $\ZZ$, respectively, and the ``unit circle''
group $\TT$ is the group $(-\frac{1}{2},\frac{1}{2}]$ with addition
mod~$1$. Except when we specify otherwise, these groups carry their
usual metrizable topology.

The symbol $\PP$ denotes the set of positive prime integers.

The set of homomorphisms $h:G\rightarrow\TT$, a group under pointwise
operation, is denoted $\Hom(G,\TT)$.
For a subgroup $A$ of $\Hom(G,\TT)$ we denote by $(G,\sT_A)$ the group
$G$ with the topology induced by $A$. Evidently $(G,\sT_A)$ is a
Hausdorff topological group if and only if $A$ separates points of $G$.
The topology $\sT_A$ is the coarsest topology on $G$ for which
the homomorphism $e_A:G\rightarrow\TT^A$ given by
$(e_A(x))_h=h(x)$ ($x\in G$, $h\in A$) is continuous. When $G=(G,\sT)$
is a topological group, the set of $\sT$-continuous functions in
$\Hom(G,\TT)$ is a subgroup of $\Hom(G,\TT)$ denoted 
$\widehat{G}$ or $\widehat{(G,\sT)}$; in this
case the topology $\sT_{\widehat{G}}$ is the {\it Bohr topology}
associated with $\sT$, and $(G,\sT_{\widehat{G}})$ is denoted $G^+$
or $(G,\sT)^+$. When $\widehat{(G,\sT)}$ separates points we say that
$G$ is a {\it maximally
almost periodic} group and we write $G=(G,\sT)\in{\bf
MAP}$. Whether or not $(G,\sT)\in{\bf MAP}$, the closure
of $e[G]$ in $\TT^{\widehat{G}}$, denoted $b(G)$ or $b(G,\sT)$, is the {\it
Bohr compactification} of $(G,\sT)$.

The Bohr compactification $b(G)$ of a topological
group $G$ is characterized by the
condition that each continuous homomorphism from $G$ into a compact
Hausdorff group extends continuously to a homomorphism from $b(G)$. From 
this and the uniform continuity of continuous homomorphisms it
follows that if $D$ is a dense subgroup of $G$ then $b(D)=b(G)$. It
is conventional to suppress mention of the function
$e_{\widehat{G}}$ and to write simply $\widehat{G}=\widehat{G^+}$. When
$G\in{\bf MAP}$ we write $G^+\subseteq b(G)\subseteq\TT^{\widehat{G}}$,
the inclusions being both algebraic and topological.

A group $G$ with its discrete topology is denoted $G_d$. For notational
convenience, and following van Douwen~\cite{vdii}, for every (Abelian)
group $G$ we write $G^\#=(G_d)^+\subseteq b(G_d)$.

A subset $S$ of a topological group $G$ is said to be
{\it bounded} in $G$ if for every nonempty open
$V\subseteq G$ there is finite $F\subseteq G$ such that $S\subseteq F+V$;
$G$ is {\it locally bounded} [resp., {\it totally bounded}] if some
nonempty open subset of $G$ is bounded [resp., $G$ itself is bounded]. It 
is a theorem
of Weil~\cite{weili} that each locally bounded group $G$ embeds as a dense
topological subgroup of a locally compact group $W(G)$, unique in the obvious
sense; the group $W(G)$ is compact if and only if $G$ is totally bounded. 
We denote by {\bf LCA} [resp., \textbf{LBA}] the class of locally
compact
[resp., locally bounded] Hausdorff Abelian groups. The relation ${\bf
LCA}\subseteq{\bf MAP}$
is a well known consequence of the Gel$'$fand-Ra\u{\i}kov Theorem
(cf.~\cite[22.17]{hri}); since each subgroup $S\subseteq G\in{\bf
MAP}$ clearly satisfies $S\in{\bf MAP}$, we have in fact the relations
${\bf LCA}\subseteq\textbf{LBA}\subseteq{\bf MAP}$.

\begin{lemma}
Let $S$ be a subgroup of $G\in$~{\em {\bf LBA}}. 
Then
\begin{itemize}
\item[(a)]
$S$ is {\it dual-embedded} in $G$ in the sense that each
$h\in\widehat{S}$ extends to an element of $\widehat{G}$;
\item[(b)]
if $h\in\widehat{S}$ and $x\in G\backslash\overline{S}^G$, the
extension $k\in\widehat{G}$ of $h$ may be chosen so that $k(x)\neq0$.
\end{itemize}
\end{lemma}

It follows that for each subgroup $S$ of a group
$G\in\textbf{LBA}$ the topology of $S^+$ coincides with the topology
inherited by $S$ from $G^+$.
This validates the following notational convention. 
For $S\subseteq G\in{\bf LBA}$, $S$ not necessarily a subgroup of $G$, we
denote
by $S^+$ the set $S$ with the topology inherited from $G^+$.
When $G$
is discrete, so that $G^+=G^\#$, we write $S^\#$ in place of $S^+$ when
$S\subseteq G$.

\begin{theorem}
(Glicksberg~\cite{glickii}). 
Let $K\subseteq G\in{\bf LBA}$.
Then $K\in\sK(G)$ if and only if $K^+\in\sK(G^+)$.
Hence if $K\in\sK(G)$, then $K$ and $K^+$ are homeomorphic.
\end{theorem}

\begin{theorem}
(Flor~\cite{flor}. See also Reid~\cite{reid}). 
Let $G\in{\bf LBA}$ and let
$x_n\rightarrow p\in b(G)=b(W(G))$ with each $x_n\in
G^+\subseteq(W(G))^+\subseteq b(G)$. Then
\begin{itemize}
\item[(a)]
$p\in (W(G))^+$, and
\item[(b)]
not only $x_n\rightarrow p$ in $(W(G))^+\subseteq b(G)$ but also
$x_n\rightarrow p$ in $W(G)$.
\end{itemize}
\end{theorem}

\begin{remark}
Strictly speaking, the papers cited above in connection with Theorems~0.2
and 0.3 deal with groups $G\in{\bf LCA}$.
Our modest generalization to the case $G\in{\bf LBA}$
is justified by 0.2 and 0.3 as originally given
and by these facts about $G\in{\bf LBA}$:
\begin{itemize}
\item[(i)]
$G$ is a (dense) topological subgroup of $W(G)\in{\bf LCA}$; 
\item[(ii)]
$G^+$ is a (dense)
topological subgroup of $(W(G))^+$; and 
\item[(iii)]
$b(G)=b(W(G))$.
\end{itemize}
\end{remark}

In what follows, groups of the form $\widehat{G}$ will be given the
{\it compact-open} topology. This is defined as usual: the family
$$\{U(K,\epsilon):K\in\sK(G), \epsilon>0\}$$ 
is a base at
$0\in\widehat{G}$, where for $A\subseteq G$ one writes
$$U(A,\epsilon)=\{h\in\widehat{G}:x\in A\Rightarrow|h(x)|<\epsilon\}.$$

We have noted already that for $G\in{\bf MAP}$ the groups $\widehat{G}$
and $\widehat{G^+}$ are identical; that is, $\widehat{G}=\widehat{G^+}$ as
groups. 
Our principal interest in Theorem~0.2 is that for $G\in{\bf LBA}$ it gives
a topological consequence, as follows.

\begin{corollary}
Let $G\in{\bf LBA}$. 
Then $\widehat{G}=\widehat{G^+}$ as topological groups. 
That is, the compact-open topology on $\widehat{G}$ determined by 
$\sK(G)$ coincides with the compact-open topology on $\widehat{G}$ determined
by $\sK(G^+)$.
\end{corollary}

\begin{ack}
We acknowledge with thanks several helpful conversations relating to \S6
with these mathematicians: Adam Fieldsteel, Michael Hru{\v{s}}{\'{a}}k,
and Stevo Todor{\v{c}}evi{\'{c}}. 
\end{ack}

\section{The Groups $\widehat{G}$ for $G$ Metrizable}

If $D$ is a dense subgroup of an Abelian
topological group $G=(G,\sT)$
then every $h\in\widehat{D}$ extends (uniquely) to an element of
$\widehat{G}$; and of course, each $h\in\widehat{G}$ satisfies
$h|D\in\widehat{D}$. Accordingly, abusing notation slightly, we have
$\widehat{D}=\widehat{G}$ as groups. 
Since groups of the form $\widehat{G}$ carry the compact-open topology, it
is natural to inquire whether the identity $\widehat{G}=\widehat{D}$ is
topological as well as algebraic. 
Informally: do $\sK(G)$ and $\sK(D)$ induce the same topology on the set
$\widehat{G}=\widehat{D}$? 
The question provokes this definition.

\begin{definition}
Let $G$ be an Abelian topological group.
\begin{itemize}
\item[(a)]
Let $D$ be a dense subgroup of $G$. Then $D$ {\it determines} $G$
(alternatively: $G$ is {\it determined} by $D$) if
$\widehat{G}=\widehat{D}$ as topological groups.
\item[(b)]
$G$ is {\it determined} if every dense subgroup of $G$
determines $G$.
\end{itemize}
\end{definition}

\begin{remarks}
\mbox{}
\begin{itemize}
\item[(a)]
It is a theorem of
Kaplan~\cite[2.9]{kaplan48} (cf.\ Banasczyzk~\cite[1.3]{banaszcz} and
Raczkowski~\cite[\S3.1]{raczphd} for alternative treatments, and
Au{\ss}enhofer \cite{Au98} and \cite[3.4]{Au99} for a generalization) that
for each $G$ the family $\{U(K,\frac{1}{4}):K\in\sK(G)\}$ is basic at
$0\in\widehat{G}$. 
(For notational simplicity, henceforth we write 
$U(K):=U(K,\frac{1}{4})\subseteq\widehat{G}$ for $K\in\sK(G)$.) 
Thus the condition that a group $G$ is determined by its dense subgroup
$D$ reduces to (i.e., is equivalent to) the condition that $\sK(D)$ is
cofinal in $\sK(G)$ in the sense that for each $K\in\sK(G)$ there is
$E\in\sK(D)$ such that $U(E)\subseteq U(K)$.
\item[(b)]
Let $D$ and $S$ be dense subgroups of a topological group $G$ such that
$D\subseteq S\subseteq G$. 
Then since $\sK(D)\subseteq\sK(S)\subseteq\sK(G)$, $D$
determines $G$ if and only if $D$ determines $S$ and $S$ determines
$G$. In particular, a dense subgroup of a determined group is determined.
\item[(c)]
The principal theorem in this corner of mathematics is the following
result, obtained independently by Au\ss enhofer \cite[4.3]{Au99} and 
Chasco~\cite{chasco}. 
This is the point of departure of the present inquiry. 
\end{itemize}
\end{remarks}

\begin{theorem}
Every metrizable, Abelian group is determined.
\end{theorem}

\begin{discussion}
Is every topological group determined? Is every {\bf MAP} group 
determined? Are there nonmetrizable, determined groups? Is every closed
(or, open) subgroup of a determined group itself determined? Is the class of
determined groups closed under passage to continuous homomorphisms?
Continuous isomorphisms? The formation of products? These are some of
the questions we address.
\end{discussion}

\section{Determined Groups: $G$ vs. $G^+$}

\begin{lemma} 
Let $D$ be a subgroup of $G\in{\bf LBA}$. 
Then $D$ is dense in $G$ if and only if $D^+$ is dense in $G^+$.
\end{lemma}

\begin{corollary}
Let $D$ be a subgroup of $G\in{\bf LBA}$. Then $D$
determines
$G$ if and only if $D^+$ determines $G^+$.
\end{corollary}

\begin{theorem}
Let $G\in{\bf LBA}$. 
Then $G$ is determined if and only if $G^+$ is determined.
\end{theorem}

\begin{theorem}
Let $G$ be an {\bf LBA} group such that $G^+$ determines $b(G)$. 
Then
\begin{itemize}
\item[(a)]
$G$ is totally bounded (and hence $G=G^+$); and
\item[(b)]
if also $G\in{\bf LCA}$ then $G$ is compact (and hence
$G=G^+=b(G)$).
\end{itemize}
\end{theorem}

\begin{corollary}
Let $G\in{\bf LBA}$. 
Then $b(G)$ is determined if and only if $W(G)$ is compact and
determined; in this case $W(G)=b(G)$.
\end{corollary}

\begin{corollary}
Let $G\in{\bf LCA}$. If $G$ is noncompact then $G^+$
does not determine $b(G)$ (and hence $b(G)$ is not determined).
\end{corollary}

\begin{theorem}
Let $G$ be a closed subgroup of a product of {\bf LBA} groups. 
Then a dense subgroup $D$ of $G$ determines $G$ if and only $D^+$ 
determines $G^+$. 
Thus $G$ is determined if and only if $G^+$ is determined.
\end{theorem}

\section{Determined Groups: Some Examples}

\begin{theorem}
There are totally bounded, nonmetrizable, determined groups.
\end{theorem}

\begin{proof}
Let $G$ be an arbitrary determined {\bf LBA} group such
that $G$ is not totally bounded. (Appealing to Theorem~1.3, one might
choose $G\in\{\ZZ,\QQ,\RR\}$.) That $G^+$ is as required follows from
three facts: 
\begin{itemize}
\item[(a)]
$G^+$ is determined~(Theorem~2.3); 
\item[(b)]
a group with a dense metrizable subgroup is itself metrizable 
\cite[Prop. IX \S2.1.1]{Bo2-66};
\item[(c)]
$b(G)$ is not metrizable.
\end{itemize}
\end{proof}

\begin{theorem}
A nondetermined group may have a dense, determined subgroup.
\end{theorem}

\begin{theorem}
The image of a nondetermined group under a continuous homomorphism may be
determined.
\end{theorem}

\begin{proof} We see in Theorem~4.8 below that compact groups of weight $\geq\cc$
are nondetermined. Each such group maps by a continuous homomorphism onto
either the group $\TT$ or a group 
of the form $(\ZZ(p))^\omega$ ($p\in\PP$) \cite{comfremusiv}, and such groups are determined 
by Theorem~1.3.
\end{proof}

\begin{discussion}
Obviously an {\bf LBA} group with no proper dense subgroup is vacuously
determined. We mention three classes of such groups.
\begin{itemize}
\item[(i)]
Discrete groups.
\item[(ii)]
Groups of the form $G^\#=(G_d)^+$. 
(It is well known \cite[2.1]{comfsaks} that every subgroup of such a group
is closed.) 
\item[(iii)]
${\bf LCA}$ groups of the type given by Rajagopalan and 
Subrahmanian \cite{rajsub}. 
Specifically, let $\kappa\geq\omega$, fix $p\in\PP$, and topologize the
group $G:=(\ZZ(p^\infty))^\kappa$ so that its subgroup 
$H:=(\ZZ(p))^\kappa$ in its usual compact topology is open-and-closed in
$G$. 
\end{itemize}
\end{discussion}

\begin{theorem}
\mbox{}
\begin{itemize}
\item[(a)]
A determined group may contain a nondetermined open-and-closed subgroup.
\item[(b)]
There are non-totally bounded, nonmetrizable, determined {\rm 
\textbf{LBA}} groups.
\end{itemize}
\end{theorem}

Although compact groups of the form $K^\kappa$ with $\kappa\geq\cc$ are
not determined, we see in Corollary~3.11 below that such groups do
contain nontrivial determining subgroups.

\begin{notation}
Let $\{G_i:i\in I\}$ be a set of groups, let $S_i\subseteq G_i$, and let
$p\in G:=\Pi_{i\in I}\,G_i$. 
Then
\begin{itemize}
\item[(i)]
$s(p)=\{i\in I:p_i\neq0_i\}$;
\item[(ii)]
$\oplus_{i\in I}\,G_i=\{x\in G:|s(x)|<\omega\}$; and
\item[(iii)]
$\oplus_{i\in I}\,S_i=(\Pi_{i\in I}\,S_i)\cap(\oplus_{i\in
I}\,G_i)$.
\end{itemize}
In this context we often identify $S_i$ with the subset
$S_i\times\{0_{I\backslash\{i\}}\}$ of $G$. In particular we
write $G_i\subseteq G$ and we identify $\widehat{G_i}$ with
$\{h|G_i:h\in\widehat{G}\}$. 
\end{notation}

We use the following property to find some determining subgroups of
certain (nondetermined) products.

\begin{definition}
A topological group $G$ has the {\it cofinally zero} property if for all
$K\in\sK(G)$ there is $F\in\sK(G)$ such that every
$h\in U(F)$ satisfies $h|K\equiv0$.
\end{definition}

\begin{remark}
We record two classes of groups with the cofinally zero property.
\begin{itemize}
\item[(i)]
$G$ is a determining subgroup of a compact Abelian group. 
(There is $F\in\sK(G)$ such that $U(F)=\{0\}$, so each $h\in U(F)$
satisfies $h|K\equiv0$ for all $K\in\sK(G)$.)
\item[(ii)]
$G$ is a torsion group of bounded order. 
(Given $K\in\sK(G)$, let $n>4$ satisfy $nx=0$ for all $x\in G$ and use
Remark 1.2(a) to choose $F\in\sK(G)$ such that 
$U(F)\subseteq U(K,\frac{1}{n})$.)
\end{itemize}
\end{remark}

\begin{lemma}
Let $\{G_i:i\in I\}$ be a set of {\bf LBA} groups with the cofinally zero
property and let $G=\Pi_{i\in I}\,G_i$. If $D_i$ is a
dense, determining subgroup of $G_i$, then $D:=\oplus_{i\in I}\,D_i$
determines $G$.
\end{lemma}

\begin{corollary}
Let $\{G_i:i\in I\}$ be a set of determined {\bf LBA} groups with the
cofinally zero property and let $G=\Pi_{i\in I}\,G_i$.
If $D_i$ is a dense subgroup of $G_i$, then $\oplus_{i\in I}\,D_i$
determines $G$.
\end{corollary}

\begin{corollary}
Let $\{G_i:i\in I\}$ be a set of compact Hausdorff groups and let
$G=\Pi_{i\in I}\,G_i$.
Then $\oplus_{i\in I}\,G_i$ determines $G$.
\end{corollary}

\begin{theorem} \label{3.12}
The image under a continuous homomorphism of a compact determined group is
determined.
\end{theorem}

\begin{remark}
It is easily checked that if a locally compact space $X$ is
$\sigma$-compact then it is {\em hemicompact}, i.e., 
some countable subfamily $\{K_n:n<\omega\}$ of $\sK(X)$ is cofinal in
$\sK(X)$ in the sense that for each $K\in\sK(X)$ there is $n<\omega$ such
that $K\subseteq K_n$. 
It follows that if an {\bf LCA} group $G$ is $\sigma$-compact
(equivalently: Lindel\"of) then $w(\widehat{G})\leq\omega$, so
$\widehat{G}$ in this case is determined by Theorem~1.3.
\end{remark}

\section{Nondetermined groups: Some Examples}

The principal result of this section is that compact Abelian groups of weight
$\geq\cc$ are nondetermined. 

\begin{lemma}
Let $G$ be an {\bf LBA} group with a proper dense subgroup $D$
such that each $K\in\sK(D)$ satisfies either
\begin{itemize}
\item[(i)]
$K$ is finite or
\item[(ii)]
$\langle K\rangle$ is closed in $D$.
\end{itemize}
Then $D$ does not determine $G$.
\end{lemma}

We noted in Theorem 2.4(b)
that if $G^+$ determines $b(G)$ with $G\in{\bf LCA}$, then $G$ is
compact (in fact $G=G^+=b(G)$).
Lemma~4.1 allows a more direct proof in the case that $G$ is discrete.

\begin{corollary}
Let $G$ be an infinite Abelian group. 
Then $G^\#$ does not determine $b(G_d)$.
\end{corollary}

\begin{lemma}
(\cite{comfrossi}). 
Let $G$ be an Abelian group.
\begin{itemize}
\item[(a)]
If $A$ is a point-separating subgroup of $\Hom(G,\TT)$, then $(G,\sT_A)$
is a totally bounded, Hausdorff topological group with $\widehat{(G,\sT_A)}=A$;
\item[(b)] 
for every totally bounded Hausdorff topological group topology $\sT$ on
$G$ the subgroup $A:=\widehat{(G,\sT)}$ of $\Hom(G,\sT)$ is
point-separating and satisfies $\sT=\sT_A$.
\end{itemize}
\end{lemma}

\begin{discussion}
It is easily checked that for each Abelian group $G$ the set $\Hom(G,\TT)$
is closed in the compact space $\TT^G$. 
Thus $\Hom(G,\TT)$, like every Hausdorff (locally) compact group, carries
a Haar measure. 
Our convention here is that Haar measure is complete, so in particular
every subset of a measurable set of measure $0$ is itself measurable (and
of measure $0$).

Concerning Haar measure $\lambda$ on a {\bf LCA} group $G$ we
appeal frequently to the
{\it Steinhaus-Weil Theorem}: {\it If $S\subseteq G$ is
$\lambda$-measurable and $\lambda(S)>0$, then the difference set
$S-S:=\{x-y:x,y\in S\}$ contains a nonempty open subset of $G$; thus
$S$, if a subgroup of $G$, is open in $G$.}
\end{discussion}

\begin{lemma}
(\cite[3.10]{comftrigwu}). 
Let $G$ be an Abelian group, let $\{x_n:n<\omega\}$ be a faithfully index
sequence in $G$, and let 
$$S:=\{h\in\Hom(G,\TT):h(x_n)\rightarrow0\in\TT\}.$$
Let $\lambda$ be the Haar measure of $\Hom(G,\TT)$. Then $S$ is
a $\lambda$-measurable subgroup of $\Hom(G,\TT)$, with $\lambda(S)=0$.
\end{lemma}

\begin{theorem}
Let $X$ be a compact Hausdorff space such that
$|X|<2^{\aleph_1}$. Then
\begin{itemize}
\item[(a)]
(\cite{hajjuh76}) 
$X$ contains a closed, countably infinite subspace; and
\item[(b)] 
$X$ contains a nontrivial convergent sequence.
\end{itemize}
\end{theorem}

\begin{theorem}
Let $G$ be an Abelian group such that $|G|<2^{\aleph_1}$ and let $A$ be a
dense subgroup of $\Hom(G,\TT)$ such that either
\begin{itemize}
\item[(i)]
$A$ is non-Haar measurable, or
\item[(ii)]
$A$ is Haar measurable, with $\lambda(A)>0$.
\end{itemize}
Then $(G,\sT_A)$ does not determine $W(G,\sT_A)$.
\end{theorem}

\begin{theorem}
Let $G$ be a compact, Abelian group such that $w(G)\geq\cc$. Then $G$ is
not determined.
\end{theorem}
\begin{proof} (Outline) Step 1. $\TT^\cc$ is not determined. [Proof.
There is a nonmeasurable subgroup $A$ of $\TT$ algebraically of the form
$\oplus_\cc\,\ZZ$. Apply Theorem~4.7(i) to the (dense) subgroup
$e_A(\ZZ)$ of $\TT^A=\TT^\cc$.]

Step 2. $F^\cc$ is not determined ($F$ a finite Abelian group, $|F|>1$).
[Proof. We have $b((\oplus_\omega\,F)_d)=F^\cc$, so Corollary 4.2
applies.]

Step 3. There is a continuous isomorphism $\phi:G\twoheadrightarrow
K^\cc$ with either $K=\TT$ or $K=F$ as in Step~2, so Theorem~\ref{3.12}
applies.
\end{proof} 

\begin{corollary}
[CH] Let $G$ be a compact Abelian group. 
Then $G$ is determined if and only if $G$ is metrizable.
\end{corollary}

\begin{corollary}
[CH] Let $\{G_i:i\in I\}$ be a set of compact Abelian groups
with each $|G_i|>1$, and let $G=\Pi_{i\in I}\,G_i$. 
Then $G$ is determined if and only if $|I|\leq\omega$ and each $G_i$ is
determined.
\end{corollary}

We close this section with an example indicating that the intersection
of dense, determining subgroups may be dense and nondetermining.

\begin{theorem}
There are dense, determining subgroups $D_i$ ($i=0,1$) of
$\TT^\cc$ such that $D_0\cap D_1$ is dense in $\TT^\cc$ and does not
determine $\TT^\cc$.
\end{theorem}

\begin{proof}
Let $Z$ be a dense, cyclic, nondetermining subgroup of $\TT^\cc$, let 
$A_i$ ($i=0,1$) be
dense, torsion subgroups of $\TT$ such that $A_0\cap A_1=\{0_\TT\}$, and
set $D_i:=Z+\oplus_\cc\,A_i\subseteq\TT^\cc$ ($i=0,1$). Then $A_i$
determines $\TT$ by Theorem~1.3 so $\oplus_\cc\,A_i$ determines
$\TT^\cc$ by Lemma~3.9, so $D_i$ determines $\TT^\cc$ ($i=0,1$); but the dense
subgroup $Z=D_0\cap D_1$ of $\TT^\cc$ does not determine $\TT^\cc$.
\end{proof}

\section{Concerning Topological Linear Spaces}

\begin{remark}
Let $\kappa$ be a cardinal number and denote by $l^1_\kappa$ the space of
real $\kappa$-sequences $x=\{x_\xi:\xi<\kappa\}$ such that 
$||x||_1:=\sum_{\xi<\kappa}|x_\xi|<\infty$.
The additive topological group $l_\kappa^1$ respects compactness
(cf.\ Remus and Trigos-Arrieta~\cite{remustrig97a}).

The group $\widehat{(l_\kappa^1)^+}$ is not discrete, so the Weil
completion $W((l_\kappa^1)^+)$ is another example of a compact
nondetermined group.
\end{remark}

\begin{mydefinition} 
A topological group $G$ is {\it (group) reflective}
if the evaluation mapping
$\Omega_G:G\rightarrow \widehat{\widehat{G}}$
defined by
$\Omega_G(x)(h):=h(x)$ for $x\in G$, $h\in\widehat{G}$
is a topological isomorphism
of $G$ onto $\widehat{\widehat G}$.
\end{mydefinition}

\begin{theorem}
Let $G$ be a noncomplete, reflective group and let $R(G)$ be its Ra\u{\i}kov
completion. 
Then $R(G)\in\bf MAP$ and $G$ does not determine $R(G)$.
\end{theorem}

Example 5.4 {\em infra} illustrates Theorem 5.2.

\begin{definition}
A reflexive locally convex vector space (LCS) in which every closed
bounded subset is compact is called a
{\it Montel space}.

Reflexivity and boundedness (Schaefer\cite{shaf:86} \S I.5, \S IV.5) are 
meant here in the sense of topological vector spaces. 
By a {\it Montel group} we mean the underlying (additive)
topological group of a Montel space. Since by definition these are reflexive
LCS, Montel groups are reflective as
proven by Smith~\cite{smith52}.
\end{definition}

\begin{example}
K\={o}mura~\cite{komura64} and Amemiya and
K\={o}mura~\cite{amemiyakomura} construct by induction three
different noncomplete Montel spaces, the completion of each being
a ``big product'' of copies of $\RR$, and one of them being exactly
$\RR^\cc$.
These groups indicate that
Theorem~5.2 is not vacuous. One of the groups constructed
in \cite{amemiyakomura} is
separable. Thus in particular, again by Theorem~5.2,
we see that $\RR^\cc$ has a countable dense subgroup which does not determine
$\RR^\cc$.
\end{example}

The remarks above show again that the property of being
determined is not $\cc$-productive.

\section{Cardinals $\kappa$ Such That $\omega<\kappa\leq\cc$}

It is well known (cf.\ for example \cite[2.18]{kunen80} or 
\cite[8.2.4]{cies}) that under Martin's Axiom [MA] every cardinal
$\kappa$ with $\omega\leq\kappa<\cc$ satisfies $2^\kappa=\cc$.
In particular under MA~$+\neg$CH
it follows from Theorem~4.6(b)
that every compact Hausdorff
space $X$ such that $|X|< 2^{\aleph_1}=\cc$ contains a nontrivial convergent
sequence.
Malykhin and {\v{S}}apirovski{\u{\i}}~\cite{malysap73} have achieved
a nontrivial extension of this result: Under MA, every compact Hausdorff
space $X$ with $|X|\leq\cc$ contains a nontrivial convergent sequence.

\begin{theorem}
[MA] Let $G$ be a group with $|G|\leq 2^{\omega}$, and let $A$ be a dense
nonmeasurable subgroup of $\widehat{G_d}$. 
Then every compact subset of $(G,\sT_A)$ is finite, so its completion
$W(G,\sT_A)$ is not determined.
\end{theorem}

If we denote by $\lambda_G$ the (completed) Haar measure on a LCA group $G$, let 
$\lambda^*_{G}$ stand for the associated outer measure. The existence of a
nonmeasurable subset $X$ of $\TT$ (with $|X|=\cc$) is well known, so the
case $\kappa=\cc$ of the following theorem generalizes the statement in
Step~1 of the proof of Theorem~4.8.

\begin{theorem}
Let $\omega < \kappa\leq\cc$.
If there is $X\in [\TT]^\kappa$ such that $\lambda_\TT^*(X)>0$, then there 
is a nonmeasurable, free Abelian subgroup $A$ of
$\TT$ algebraically of the form 
$A=\oplus_\kappa\,\ZZ$.
\end{theorem}

Responding to a question on a closely related matter, Stevo 
Todor\-\v{c}evi\'{c}~\cite{To01} proposed and proved the above result for 
$\kappa=\aleph_1$. 
His proof additionally yields that
$X\backslash\tor(\TT)$ can be broken into $\omega$-many pairwise disjoint
independent sets, each of cardinality $\aleph_1$.

For torsion groups of prime order, we obtain the following.

\begin{theorem}
Let $F$ be a finite group of prime order $p$, and let $\kappa_1,
\kappa_2$ be infinite cardinals such that 
$\kappa_1\leq2^{\kappa_2}$. 
If there is $X\in [F^{\kappa_2}]^{\kappa_1}$ such that $\lambda^*_{F^{\kappa_2}}(X)>0$,
then there is a nonmeasurable subgroup $A$ of $F^{\kappa_2}$ algebraically
of the form $A=\oplus_{\kappa_1} F$.
\end{theorem}

\begin{discussion}
For an ideal $\sI$ of subsets of a set $S$ we write as usual
$${\mathrm{non}}(\sI)=\min\{|Y|:Y\subseteq S,~Y\notin\sI\}.$$
Let $F$ be a finite group ($|F|>1$),
let $\sN(\TT)$ and $\sN(F^\omega)$ denote
the $\sigma$-algebra of $\lambda_\TT$- and
$\lambda_{F^\omega}$-measurable sets of measure zero. As with any two
compact metric spaces of equal cardinality equipped with atomless
(``continuous'') probability measures, the spaces $\TT$ and $F^\omega$
are Borel-isomorphic in the sense that there is a bijection
$\phi:\TT\twoheadrightarrow F^\omega$ 
such that the associated bijection
$\overline{\phi}:\sP(\TT)\twoheadrightarrow\sP(F^\omega)$ carries
the Borel algebra $\sB(\TT)$ onto the Borel algebra $\sB(F^\omega)$ 
in such a way that 
$\lambda_{F^\omega}(\overline{\phi}(B))=\lambda_\TT(B)$ for each $B\in\sB(\TT)$.
(See \cite[17.41]{kechris} or \cite[3.4.23]{sriva98} for a proof of
this ``Borel isomorphism Theorem for measures''.) 
\end{discussion}

\begin{lemma}
The cardinals $\non(\sN(\TT))$ and $\non(\sN(F^\omega))$ are equal.
\end{lemma}

We write
$$\non(\sN):=\non(\sN(\TT))=\non(\sN(F^\omega)),$$ 
a definition justified by Lemma~6.5.

\begin{theorem}
Let $G$ be a compact Abelian group such that $w(G)\geq\non(\sN)$. 
Then $G$ is nondetermined.
\end{theorem}

\section{Questions} 

\begin{question}
Is there a compact group $G$ with a countable dense subgroup $D$ such that
$w(G)>\omega$ and $D$ determines $G$?
\end{question}

\begin{question}
If $\{G_i:i\in I\}$ is a set of topological Abelian groups and $D_i$ is a
dense determining subgroup of $G_i$, must
$\oplus_{i\in I}\,D_i$ determine $\Pi_{i\in I}\,G_i$? In particular,
does $\oplus_{i\in I}\,G_i$ determine $\Pi_{i\in I}\,G_i$? In
particular, does $\oplus_\cc\,\RR$ determine $\RR^\cc$?
\end{question}

\begin{discussion}
Consider the following cardinals:
\begin{itemize}
\item[(a)]
$\mm_\TT :=$ the least cardinal $\kappa$ such that $\TT^\kappa$
is nondetermined;
\item[(b)]
$\mm_{f\exists}$ [resp., $\mm_{f\forall}] :=$ the
least cardinal $\kappa$ such that some [resp., each] finite group
$F$ has $F^\kappa$ nondetermined;
\item[(c)]
$\mm_{c\exists}$ [resp., $\mm_{c\forall}] :=$ the least
cardinal $\kappa$ such that some [resp., each] compact
abelian group of weight $\kappa$ is nondetermined;
\item[(d)]
$\mm_{p\exists}$ [resp., $\mm_{p\forall}] :=$ the least cardinal
$\kappa$ such that some [resp., each] product of $\kappa$-many compact
determined groups is nondetermined.
\end{itemize}

It follows from Theorems 1.3 and 6.6 that each $\mm_x$, with the
possible exception of $\mm_{p\exists}$, satisfies
$\aleph_1 \leq \mm_x\leq\non(\sN)$.
Further if $\non(\sN)=\aleph_1$, then all seven cardinals $\mm_x$ are 
equal to $\aleph_1$. The condition $\non(\sN)=\aleph_1$ is clearly
consistent with CH, and it has been shown to be consistent as well
with $\neg$CH (see for example \cite{bartjudah}, \cite{fremlin84} and
\cite[Example 1, page 568]{jechii}), so in particular there are models of
ZFC~+~$\neg$CH in which every compact (Abelian) group $G$ satisfies: $G$
is determined if and only if $G$ is metrizable. (Without appealing
to the cardinal $\non(\sN)$,
Michael Hru{\v{s}}{\'{a}}k~\cite{Hr00} in informal conversation
suggested the existence of models of ZFC~+~$\neg$CH in which
$\{0,1\}^{\aleph_1}$ is nondetermined.)
\end{discussion}

\begin{question}
Are the various cardinal numbers $\mm_x$ equal in ZFC?
Are they equal to one of the familiar ``small cardinals'' conventionally
noted in
the Cicho{\'n} diagram (cf.~\cite{bartjudah}, \cite{vaughan90})?
Is each $\mm_x=\non(\sN)$? Is each $\mm_x=\aleph_1$?
Is each $\cf(\mm_x)>\omega$? 
\end{question}

We know of no models of ZFC in which $\TT^{\aleph_1}$,
or some group of the form $F^{\aleph_1}$ ($F$ finite, $|F|>1$), is
determined, so we are forced to consider the
possibility that the following questions have an affirmative answer.

\begin{question}
Are the following (equivalent) statements theorems of ZFC?
\begin{itemize}
\item[(a)]
The group $\TT^{\aleph_1}$ and groups of the form $F^{\aleph_1}$ ($F$
finite, $|F|>1$) are nondetermined.
\item[(b)]
A compact abelian group $G$ is determined if and only if $G$ is 
metrizable. 
\end{itemize}
\end{question}

The following question is suggested by those above.

\begin{question}
Is there in ZFC a cardinal $\kappa$ such that a compact group $G$ is
determined if and only if $w(G)<\kappa$?
\end{question}

\begin{question}
Is it consistent with ZFC that $\mm_{p\exists}=2$?
Is it consistent with ZFC that $\mm_{p\exists}=\omega$?
\end{question}

Question 7.7 has analogues in the context of groups which are not
assumed to be compact, as follows.

\begin{question}
In ZFC alone or in augmented axiom systems: Is the product of finitely
many determined groups necessarily determined? If $G$ is determined, is
$G\times G$ necessarily determined?
\end{question}

%\bibliographystyle{amsplain}
%\bibliography{04}
\providecommand{\bysame}{\leavevmode\hbox to3em{\hrulefill}\thinspace}
\providecommand{\MR}{\relax\ifhmode\unskip\space\fi MR }
\providecommand{\MRhref}[2]{%
  \href{http://www.ams.org/mathscinet-getitem?mr=#1}{#2}
}
\providecommand{\href}[2]{#2}

\end{document}